\documentclass{article}
\usepackage{amssymb}
\usepackage{graphicx}
\usepackage{amsmath}
\usepackage{color}

\input{tcilatex}
\begin{document}

\title{Plancherel Theorem and the Left Ideals of$\ $the Group Algebra for
the Jacobi Group.}
\author{Kahar El-Hussein \\
Department of Mathematics, Faculty of Science, \\
Al Furat University, Dear El Zore, Syria \ \ and\\
Department of Mathematics, Faculty of Arts Science Al Quryyat, \\
Al-Jouf University, KSA \\
\textit{E-mail : kumath@ju.edu.sa, kumath@hotmail.com }}
\maketitle

\begin{abstract}
Let $G=SL(2,\mathbb{R})$\ be the $2\times 2$ connected real semisimple Lie
group and let $KAN$ be the Iwasawa decomposition of $SL(2,\mathbb{R})$. Let $%
J=H\rtimes SL(2,$ $\mathbb{R})$ be the Jacobi group, which is the semidirect
product of the two groups $H$ with $SL(2,$ $\mathbb{R})$. It plays an
important role in Quantum Mechanics. The purpose of this paper is to define
the Fourier transform in order to obtain the Plancherel theorem for the
group $J$. To this end a classification of all left ideals of the group
algebra $L^{1}(H\rtimes AN)$.
\end{abstract}

\bigskip \textbf{Keywords}: Jacobi Group, Iwasawa Decomposition, Fourier
Transform and Plancherel Theorem, Left Ideals

\textbf{AMS 2000 Subject Classification:} $43A30\&35D$ $05$

\section{\textbf{\ Introduction}}

\bigskip \textbf{1.1. }The Jacobi group the semidirect product of the
Heisenberg and the symplectic group $SL(2,\mathbb{R})$ plays an important
role in quantum mechanics. In Quantum optics represent a physical
realization of the coherent states associated to the Jacobi group. The
Jacobi group is responsible for the squeezed states and has an important
object in quantum mechanics, geometric quantization, optics. Abstract
harmonic analysis is the field of the most modern branches of harmonic
analysis, having its roots in the mid-twentieth century, is analysis on
topological groups. If the group is neither abelian nor compact, no general
satisfactory theory is currently known.

\textbf{1.2.} First In this paper I will define the Fourier transform in
order to establish the Plancherel formula on the Jacobi group $H\rtimes SL(2,%
\mathbb{R}),$ where $H$ is the $3$-dimensional Heisenberg group and 
\begin{equation}
SL(2,\mathbb{R})=\{\left( X=%
\begin{array}{cc}
a & b \\ 
c & d%
\end{array}%
\right) :\text{ }detX=1\}
\end{equation}

Secondly, I will give classification for all left ideals of the group
algebra $L^{1}(H\rtimes AN),$ where\ $AN$ is the solvable Lie group in the
Iwasawa decomposition $KAN$ of $SL(2,\mathbb{R}).$

\section{Fourier Transform and Plancherel Formula on $SL(2,\mathbb{R})$}

\textbf{2.1. }In the following and far away from the representations theory
of Lie groups we use the Iwasawa decomposition of $SL(2,\mathbb{R}),$ to
define the Fourier transform and to demonstrate Plancherel formula on the
connected real semisimple Lie group $SL(2,\mathbb{R}).$ Therefore let $SL(2,%
\mathbb{R})$\ be the real Lie group, which is

\begin{equation}
SL(2,\mathbb{R})=\{\left( 
\begin{array}{cc}
a & b \\ 
c & d%
\end{array}%
\right) :\text{ }(a,b,c,d)\in \mathbb{R}^{4}\text{ \textit{and} }ad-bc=1\}
\end{equation}%
and let $SL(2,\mathbb{R})=KNA$ be the Iwasawa decomposition of $SL(2,\mathbb{%
R})$, where 
\begin{eqnarray}
K &=&\{\left( 
\begin{array}{cc}
\cos \phi & -\sin \phi \\ 
\sin \phi & \cos \phi%
\end{array}%
\right) =SO(2):\text{ }\phi \in \mathbb{R}\text{ \ }\}  \notag \\
N &=&\{\left( 
\begin{array}{cc}
1 & n \\ 
0 & 1%
\end{array}%
\right) :\text{ }n\in \mathbb{R}\text{ }\}  \notag \\
A &=&\{\left( 
\begin{array}{cc}
a & 0 \\ 
0 & a^{-1}%
\end{array}%
\right) :\text{ }a\in \mathbb{R}_{+}^{\star }\}
\end{eqnarray}

\bigskip Hence every $g\in SL(2,\mathbb{R})$ can be written as $g=kan\in
SL(2,\mathbb{R}),$ where $k\in K,$ $a\in A,$ $n\in \mathbb{R}.$

\textbf{2.2.} We denote by $L^{1}(G)$ the Banach algebra that consists of
all complex valued functions on the group $G$, which are integrable with
respect to the Haar measure $dg$ of $G$ and multiplication is defined by
convolution product on $G$ , where $G=SL(2,\mathbb{R}).$ And denote by $%
L^{2}(G)$ the Hilbert space of $G$. So we have for any $f\in L^{1}(G)$ and $%
\phi \in L^{1}(G)$ 
\begin{equation}
\phi \ast f(h)=\int\limits_{G}f(g^{-1}h)\phi (g)dg
\end{equation}%
The Haar measure $dg$ on a connected real semi-simple Lie group $G$ $=SL(n,%
\mathbb{R)}$, can be calculated from the Haar measures $dn,$ $da$ and $dk$
on $N;A$ and $K;$respectively, by the formula%
\begin{equation}
\int\limits_{G}f(g)dg=\int\limits_{A}\int\limits_{N}\int%
\limits_{K}f(ank)dadndk
\end{equation}

Keeping in mind that $a^{-2\rho }$ is the modulus of the automorphism $%
n\rightarrow $ $ana^{-1}$ of $N$ we get also the following representation of 
$dg$ 
\begin{equation}
\int\limits_{G}f(g)dg=\int\limits_{A}\int\limits_{N}\int%
\limits_{K}f(ank)dadndk=\int\limits_{N}\int\limits_{A}\int%
\limits_{K}f(nak)a^{-2\rho }dndadk
\end{equation}%
where%
\begin{equation*}
\rho =2^{-1}\sum_{\alpha \geq 0,\alpha \neq 0}m(\alpha )\alpha
\end{equation*}%
and $m(\alpha )$ denotes the multiplicity of the root $\alpha $ see $[17]$
or again $\rho =$ the dimension of the nilpotent group $N.$ Furthermore,
using the relation $\int\limits_{G}f(g)dg=\int\limits_{G}f(g^{-1})dg,$ we
receive 
\begin{equation}
\int\limits_{G}f(g)dg=\int\limits_{K}\int\limits_{A}\int%
\limits_{N}f(kan)a^{2\rho }dndadk
\end{equation}

\textbf{2.3. }Let $\Gamma $ be a connected compact Lie group and let $%
\underline{k}$ be the Lie algebra of $\Gamma $. Let $%
(X_{1},X_{2},.....,X_{m})$ a basis of $\underline{k}$ , such that the both
operators%
\begin{equation}
\Delta =\sum\limits_{i=1}^{m}X_{i}^{2}
\end{equation}%
\begin{equation}
D_{q}=\sum\limits_{0\leq l\leq q}\left(
-\sum\limits_{i=1}^{m}X_{i}^{2}\right) ^{l}
\end{equation}%
are left and right invariant (bi-invariant) on $\Gamma ,$ this basis exist
see $[2,$ $p.564)$. For $l\in 
\mathbb{N}
$, let $D^{l}=(1-\Delta )^{l}$, then the family of semi-norms $\{\sigma _{l}$%
, $l\in 
\mathbb{N}
\}$ such that%
\begin{equation}
\sigma _{l}(f)=\int_{\Gamma }\left\vert D^{l}f(y)\right\vert ^{2}dy)^{\frac{1%
}{2}},\text{ \ \ \ \ \ \ \ \ }f\in C^{\infty }(\Gamma )
\end{equation}%
define on $C^{\infty }(\Gamma )$ the same topology of the Frechet topology
defined by the semi-normas $\left\Vert X^{\alpha }f\right\Vert _{2}$ defined
as%
\begin{equation}
\left\Vert X^{\alpha }f\right\Vert _{2}=\int_{\Gamma }(\left\vert X^{\alpha
}f(y)\right\vert ^{2}dy)^{\frac{1}{2}},\text{ \ \ \ \ \ \ \ \ }f\in
C^{\infty }(\Gamma )
\end{equation}%
where $\alpha =(\alpha _{1},$.....,$\alpha _{m})\in 
\mathbb{N}
^{m},$ see $[2,p.565]$

Let $\widehat{\Gamma }$ be the set of all equivilance classes of irreducible
unitary representations of $\Gamma .$ If $\gamma \in \widehat{\Gamma }$, we
denote by $E_{\gamma }$ the space of representation $\gamma $ and $d_{\gamma
}$\ its dimension then we get\qquad

\textbf{Definition 2.1.} \textit{The Fourier transform of a function }$f\in
C^{\infty }(\Gamma )$\textit{\ is defined as} 
\begin{equation}
Tf(\gamma )=\int\limits_{\Gamma }f(x)\gamma (x^{-1})dx
\end{equation}%
\textit{where }$T$\textit{\ is the Fourier transform on} $\Gamma $

\textbf{Theorem (A. Cerezo) 2.1.} \textit{Let} $f\in C^{\infty }(\Gamma ),$ 
\textit{then we have the inversion of the Fourier transform} 
\begin{equation}
f(x)=\sum\limits_{\gamma \in \widehat{\Gamma }}d\gamma tr[Tf(\gamma )\gamma
(x)
\end{equation}

\begin{equation}
f(I_{\Gamma })=\sum\limits_{\gamma \in \widehat{\Gamma }}d\gamma
tr[Tf(\gamma )]
\end{equation}%
\textit{and the Plancherel formula} 
\begin{equation}
\left\Vert f(x)\right\Vert _{2}^{2}=\int_{\Gamma }\left\vert f(x)\right\vert
^{2}dx=\sum\limits_{\gamma \in \widehat{\Gamma }}d_{\gamma }\left\Vert
Tf(\gamma )\right\Vert _{H.S}^{2}
\end{equation}%
\textit{for any }$f\in L^{1}(\Gamma ),$ \textit{where }$I_{\Gamma }$ \textit{%
is the identity element of \ }$\Gamma $ \textit{and} $\left\Vert Tf(\gamma
)\right\Vert _{H.S}^{2}$ \textit{is the Hilbert- Schmidt norm of the operator%
} $Tf(\gamma )$\newline
Fourier did not actually assume any underlying group structure or
representation theory but we typically associate his work with the case of
the circle group in the following form using complex exponentials%
\begin{equation}
f(x)=\sum\limits_{n=-\infty }^{\infty }Tf(m)e^{ixm}=\sum\limits_{m=-\infty
}^{\infty }c_{n}e^{ixm},\text{ \ \ }m\in 
\mathbb{Z}%
\end{equation}%
where%
\begin{equation}
c_{m}=Tf(m)=\int\limits_{SO(2)}f(x)e^{-ixm}dx
\end{equation}

The group is $SO(2)=S^{1}$or $\mathbb{R}/%
\mathbb{Z}
$ and the multiplicative characters are $e^{ixn}$, group homomorphisms from
the circle $K=SO(2)$ to the multiplicative group of non-zero complex
numbers. Fourier actually preferred to express the coefficients using what
is now known as the Plancherel formula

\begin{equation}
\left\Vert f(x)\right\Vert _{2}^{2}=\int_{SO(2)}\left\vert f(x)\right\vert
^{2}dx=\sum\limits_{n=-\infty }^{\infty }\left\vert c_{m}\right\vert
^{2}=\sum\limits_{n=-\infty }^{\infty }\left\vert Tf(m)\right\vert ^{2}
\end{equation}%
where%
\begin{equation}
S^{1}=SO(2)=\{\left( 
\begin{array}{cc}
\cos \phi & -\sin \phi \\ 
\sin \phi & \cos \phi%
\end{array}%
\right) :\phi \in \mathbb{R}\ \}
\end{equation}%
\ 

\bigskip \textbf{Definition 2.2}. \textit{For any function} $f\in \mathcal{D}%
(G),$ \textit{we can define a function} $\Upsilon (f)$\textit{on }$G\times K$
$=G\times SO(2)$ \textit{by} 
\begin{equation}
\Upsilon (f)(g,k_{1})=\Upsilon (f)(kna,k_{1})=f(gk_{1})=f(knak_{1})
\end{equation}%
\textit{for }$g=kna\in G,$ \textit{and} $k_{1}\in K$ . \textit{The
restriction of} $\ \Upsilon (f)\ast \psi (g,k_{1})$ \textit{on} $K(G)$ 
\textit{is }$\Upsilon (f)\ast \psi (g,k_{1})\downarrow
_{K(G)}=f(nak_{1})=f(g)\in \mathcal{D}(G),$ \textit{and }$\Upsilon
(f)(g,k_{1})\downarrow _{K}=f(kna)$ $\in \mathcal{D}(G)$

\textbf{Remark 2.1}. $\Upsilon (f)$ is invariant in the following sense%
\begin{equation}
\Upsilon (f)(gh,h^{-1}k_{1})=\Upsilon (f)(g,k_{1})
\end{equation}

\textbf{Definition 2.3}.\textit{\ If }$f$ \textit{and }$\psi $ \textit{are
two functions belong to} $\mathcal{D}(G),$ \textit{then we can define the
convolution of } $\Upsilon (f)\ $\textit{and} $\psi $\ \textit{on} $G$ $%
\times K=G\times S^{1}=G\times SO(2)$ \textit{as}

\begin{eqnarray*}
\Upsilon (f)\ast \psi (g,k_{1}) &=&\int\limits_{G}\Upsilon
(f)(gg_{2}^{-1},k_{1})\psi (g_{2})dg_{2} \\
&=&\int\limits_{SO(2)}\int\limits_{N}\int\limits_{A}\Upsilon
(f)(knaa_{2}^{-1}n_{2}^{-1}k^{-1}k_{1})\psi
(k_{2}n_{2}a_{2})dk_{2}dn_{2}da_{2}
\end{eqnarray*}

So we get 
\begin{eqnarray*}
\Upsilon (f)\ast \psi (g,k_{1}) &\downarrow &_{K(G)}=\Upsilon (f)\ast \psi
(I_{K}na,k_{1}) \\
&=&\int\limits_{SO(2)}\int\limits_{N}\int%
\limits_{A}f(naa_{2}^{-1}n_{2}^{-1}k^{-1}k_{1})\psi
(k_{2}n_{2}a_{2})dk_{2}dn_{2}da_{2} \\
&=&\Upsilon (f)\ast \psi (na,k_{1})
\end{eqnarray*}%
where $g_{2}=k_{2}n_{2}a_{2}$

\bigskip \textbf{Definition} \textbf{2.4.} \textit{If }$f\in \mathcal{D}(G)$ 
\textit{and let} $\Upsilon (f)$ \textit{be the associated function to} $f$ , 
\textit{we define the Fourier transform of \ }$\Upsilon (f)(g,k_{1})$ 
\textit{by }%
\begin{eqnarray}
&&\mathcal{F}\Upsilon (f))(I_{S^{1}},\xi ,\lambda ,\gamma ,I_{S^{1}})=%
\mathcal{F}\Upsilon (f)(I_{S^{1}},\xi ,\lambda ,I_{S^{1}})  \notag \\
&=&\int_{S^{1}}\int_{N}\int_{A}[\sum\limits_{l=-\infty }^{\infty
}\int_{S^{1}}T\Upsilon (f)(kna,k_{1})e^{-ilk}dk]a^{-i\lambda }e^{-\text{ }%
i\langle \text{ }\xi ,\text{ }n\rangle }\text{ }e^{-imk_{1}}dadndk_{1} 
\notag \\
&=&\int_{S^{1}}\int_{N}\int_{A}[\Upsilon (f)(I_{S^{1}}na,k_{1})]a^{-i\lambda
}e^{-\text{ }i\langle \text{ }\xi ,\text{ }n\rangle }\text{ }%
e^{-imk_{1}}dadndk_{1}
\end{eqnarray}%
\textit{where }$\mathcal{F}$ \textit{is the Fourier transform on }$AN$ 
\textit{and }$T$ \textit{is the Fourier transform on} $SO(2),$ \textit{and }$%
I_{S^{1}}$ \textit{is the identity element of }$S^{1}=SO(2)$

\textbf{\textit{Plancherel's \textbf{Theorem} on the Group G \textbf{2.2.} }}%
\textit{For any function\ }$f\in $\textit{\ }$L^{1}(G)\cap $\textit{\ }$%
L^{2}(G),$\textit{we get }%
\begin{equation}
\int_{G}\left\vert f(g)\right\vert
^{2}dg=\int\limits_{A}\int\limits_{N}\int\limits_{S^{1}}\left\vert
f(kna)\right\vert ^{2}dadndk=\sum\limits_{m=-\infty }^{\infty }\int\limits_{%
\mathbb{R}}\int\limits_{\mathbb{R}}\left\Vert T\mathcal{F}f(\lambda ,\xi
,m)\right\Vert _{2}^{2}d\lambda d\xi
\end{equation}%
\textit{\ }%
\begin{equation}
f(I_{A}I_{N}I_{S^{1}})=\int\limits_{N}\int\limits_{A}\sum\limits_{m=-\infty
}^{\infty }T\mathcal{F}f((\lambda ,\xi ,m)]d\lambda d\xi
=\sum\limits_{m=-\infty }^{\infty }\int\limits_{\mathbb{R}}\int\limits_{%
\mathbb{R}}T\mathcal{F}f(\lambda ,\xi ,m)d\lambda d\xi
\end{equation}%
\textit{where} $I_{A},I_{N},$ and $I_{K}$ \textit{are the identity elements
of} $A$, $N$ \textit{and }$K$ \textit{respectively, where }$\mathcal{F}$ 
\textit{is the Fourier transform on }$AN$ \textit{and }$T$ \textit{is the
Fourier transform on} $K,$ \textit{and }$I_{K}$ \textit{is the identity
element of }$K$

\bigskip \textit{Proof: }First let $\overset{\vee }{f}$ be the function
defined by 
\begin{equation}
\ \overset{\vee }{f}(kna)=\overline{f((kna)^{-1})}=\overline{%
f(a^{-1}n^{-1}k^{-1})}
\end{equation}

Then we have%
\begin{eqnarray}
&&\int_{G}\left\vert f(g)\right\vert ^{2}dg  \notag \\
&=&\Upsilon (f)\ast \overset{\vee }{f}(I_{S^{1}}I_{N}I_{A},I_{S^{1}})  \notag
\\
&=&\int\limits_{G}\Upsilon (f)(I_{S^{1}}I_{N}I_{A}(g_{2}^{-1}),I_{S^{1}})%
\overset{\vee }{f}(g_{2})dg_{2}  \notag \\
&=&\int\limits_{A}\int\limits_{N}\int_{S^{1}}\Upsilon
(f)(a_{2}^{-1}n_{2}^{-1}k_{2}^{-1},I_{S^{1}})\overset{\vee }{f}%
(k_{2}n_{2}a_{2})da_{2}dn_{2}dk_{2}  \notag \\
&=&\int\limits_{A}\int\limits_{N}%
\int_{S^{1}}f(a_{2}^{-1}n_{2}^{-1}k_{2}^{-1})\overline{%
f((k_{2}n_{2}a_{2})^{-1})}da_{2}dn_{2}dk_{2}  \notag \\
&=&\int\limits_{A}\int\limits_{N}\int_{S^{1}}\left\vert
f(a_{2}n_{2}k_{2})\right\vert ^{2}da_{2}dn_{2}dk_{2}
\end{eqnarray}

\bigskip Secondly%
\begin{eqnarray*}
&&\Upsilon (f)\ast \overset{\vee }{f}(I_{S^{1}}I_{N}I_{A},I_{S^{1}}) \\
&=&\int\limits_{\mathbb{R}}\int\limits_{\mathbb{R}}\text{ }\mathcal{F}%
(\Upsilon (f)\ast \overset{\vee }{f})(I_{S^{1}},\lambda ,\xi
,I_{S^{1}})d\lambda d\xi \\
&=&\int_{S^{1}}\int\limits_{\mathbb{R}}\int\limits_{\mathbb{R}%
}\int\limits_{A}\int\limits_{N}\text{ }\sum\limits_{m=-\infty }^{\infty
}\sum\limits_{l=-\infty }^{\infty }\int_{S^{1}}\Upsilon (f)\ast \overset{%
\vee }{f}(kna,k_{1})e^{-ilk}dka^{-i\lambda }e^{-\text{ }i\langle \text{ }\xi
,\text{ }n\rangle }\text{ }e^{-imk_{1}}dadndk_{1}d\lambda d\xi \\
&=&\sum\limits_{m=-\infty }^{\infty }\int_{S^{1}}\int\limits_{\mathbb{R}%
}\int\limits_{\mathbb{R}}\int\limits_{A}\int\limits_{N}\Upsilon (f)\ast 
\overset{\vee }{f}(I_{S^{1}}na,k_{1})e^{-ilk}dka^{-i\lambda }e^{-\text{ }%
i\langle \text{ }\xi ,\text{ }n\rangle }\text{ }e^{-imk_{1}}dadndk_{1}d%
\lambda d\xi \\
&=&\int\limits_{\mathbb{R}}\int\limits_{\mathbb{R}}\int\limits_{A}\int%
\limits_{N}\int_{S^{1}}\int\limits_{A}\int\limits_{N}\sum\limits_{m=-\infty
}^{\infty }\int_{S^{1}}\Upsilon
(f)(I_{S^{1}}naa_{2}^{-1}n_{2}^{-1}k_{2}^{-1},k_{1})\overset{\vee }{f}%
(k_{2}n_{2}a_{2})e^{-imk_{1}}dk_{1} \\
&&dndadk_{2}dn_{2}da_{2}a^{-i\lambda }e^{-\text{ }i\langle \text{ }\xi ,%
\text{ }n\rangle }d\lambda d\xi \\
&=&\int\limits_{\mathbb{R}}\int\limits_{\mathbb{R}}\int\limits_{A}\int%
\limits_{N}\int_{S^{1}}\int\limits_{A}\int\limits_{N}\sum\limits_{m=-\infty
}^{\infty }\int_{S^{1}}f(naa_{2}^{-1}n_{2}^{-1}k_{2}^{-1}k_{1})\overset{\vee 
}{f}(k_{2}n_{2}a_{2})e^{-imk_{1}}dk_{1}dk_{2} \\
&&a^{-i\lambda }e^{-\text{ }i\langle \text{ }\xi ,\text{ }n\rangle
}dndadn_{2}da_{2}d\lambda d\xi
\end{eqnarray*}%
where%
\begin{equation}
e^{-\text{ }i\langle \text{ }\xi ,\text{ }n\rangle }=e^{-\text{ }i\text{ }%
\xi n}
\end{equation}

\bigskip\ Using the fact that%
\begin{equation}
\int\limits_{A}\int\limits_{N}\int_{S^{1}}f(kna)dadndk=\int\limits_{N}\int%
\limits_{A}\int_{S^{1}}f(kan)a^{2}dndadk
\end{equation}%
and 
\begin{eqnarray}
&&\int\limits_{\mathbb{R}}\int\limits_{A}\int\limits_{N}%
\int_{S^{1}}f(kna)e^{-\text{ }i\langle \text{ }\xi ,\text{ }n\rangle
}dadndkdd\xi  \notag \\
&=&\int\limits_{\mathbb{R}}\int\limits_{A}\int\limits_{N}%
\int_{S^{1}}f(kan)e^{-\text{ }i\langle \text{ }\xi ,\text{ }an_{1}a^{-1}%
\text{ }\rangle }a^{2}dadndkd\xi  \notag \\
&=&\int\limits_{\mathbb{R}}\int\limits_{A}\int\limits_{N}%
\int_{S^{1}}f(kan)e^{-\text{ }i\langle \text{ }a\xi a^{-1},\text{ }n\rangle
}a^{2}dadndkd\xi  \notag \\
&=&\int\limits_{\mathbb{R}}\int\limits_{A}\int\limits_{N}%
\int_{S^{1}}f(kan)e^{-\text{ }i\langle \text{ }\xi ,\text{ }n\rangle
}dadndkd\xi
\end{eqnarray}

Then we get 
\begin{eqnarray*}
&&\Upsilon (f)\ast \overset{\vee }{f}(I_{S^{1}}I_{N}I_{A},I_{S^{1}}) \\
&=&\int\limits_{\mathbb{R}}\int\limits_{\mathbb{R}}\int\limits_{A}\int%
\limits_{N}\int_{S^{1}}\int\limits_{A}\int\limits_{N}\sum\limits_{m=-\infty
}^{\infty }\int_{S^{1}}f(naa_{2}^{-1}n_{2}^{-1}k_{2}^{-1},k_{1})\overset{%
\vee }{f}(k_{2}n_{2}a_{2})e^{-imk_{1}}dk_{1}dk_{2} \\
&&a^{-i\lambda }e^{-\text{ }i\langle \text{ }\xi ,\text{ }n\rangle
}dndadn_{2}da_{2}d\lambda d\xi \\
&=&\int\limits_{\mathbb{R}}\int\limits_{\mathbb{R}}\int\limits_{A}\int%
\limits_{N}\int\limits_{A}\int\limits_{N}\sum\limits_{m=-\infty }^{\infty
}\int_{S^{1}}\int_{S^{1}}f(aa_{2}^{-1}nn_{2}^{-1}k_{2}^{-1},k_{1})\overset{%
\vee }{f}(k_{2}n_{2}a_{2})e^{-imk_{1}}dk_{1}dk_{2} \\
&&a^{-i\lambda }e^{-\text{ }i\langle \text{ }\xi ,\text{ }n\rangle
}dndadn_{2}da_{2}d\lambda d\xi \\
&=&\int\limits_{\mathbb{R}}\int\limits_{\mathbb{R}}\int\limits_{A}\int%
\limits_{N}\int\limits_{A}\int\limits_{N}\sum\limits_{m=-\infty }^{\infty
}\int_{S^{1}}\int_{S^{1}}f(ank_{2}^{-1},k_{1})\overset{\vee }{f}%
(k_{2}n_{2}a_{2})e^{-imk_{1}}dk_{1}dk_{2} \\
&&a^{-i\lambda }e^{-\text{ }i\langle \text{ }\xi ,\text{ }n\rangle
}dndadn_{2}da_{2}d\lambda d\xi \\
&=&\int\limits_{\mathbb{R}}\int\limits_{\mathbb{R}}\int\limits_{A}\int%
\limits_{N}\int\limits_{A}\int\limits_{N}\sum\limits_{m=-\infty }^{\infty
}\int_{S^{1}}\int_{S^{1}}f(ank_{2}^{-1}k_{1})\overset{\vee }{f}%
(k_{2}n_{2}a_{2})e^{-imk_{1}}dk_{1}dk_{2} \\
&&a^{-i\lambda }e^{-\text{ }i\langle \text{ }\xi ,\text{ }n\rangle
}dndadn_{2}da_{2}d\lambda d\xi \\
&=&\int\limits_{\mathbb{R}}\int\limits_{\mathbb{R}}\int\limits_{A}\int%
\limits_{N}\int\limits_{A}\int\limits_{N}\sum\limits_{m=-\infty }^{\infty
}\int_{S^{1}}\int_{S^{1}}f(ank_{1}^{-1})\overset{\vee }{f}%
(k_{2}n_{2}a_{2})e^{-imk_{1}}e^{-imk_{2}}dk_{1}dk_{2} \\
&&a^{-i\lambda }a_{2}^{-i\lambda }e^{-\text{ }i\langle \text{ }\xi ,\text{ }%
n+n_{2}\rangle }dndadn_{2}da_{2}d\lambda d\xi \\
&=&\int\limits_{\mathbb{R}}\int\limits_{\mathbb{R}}\int\limits_{A}\int%
\limits_{N}\int\limits_{A}\int\limits_{N}\sum\limits_{m=-\infty }^{\infty
}\int_{S^{1}}\int_{S^{1}}f(ank_{1}^{-1})\overline{%
f(a_{2}{}^{-1}n_{2}^{-1}k_{2}^{-1})}e^{-imk_{1}}e^{-imk_{2}}dk_{1}dk_{2} \\
&&a^{-i\lambda }e^{-\text{ }i\langle \text{ }\xi ,\text{ }n\rangle
}a_{2}^{-i\lambda }e^{-\text{ }i\langle \text{ }\xi ,\text{ }n_{2}\rangle
}dndadn_{2}da_{2}d\lambda d\xi \\
&=&\int\limits_{\mathbb{R}}\int\limits_{\mathbb{R}}\int\limits_{A}\int%
\limits_{N}\int\limits_{A}\int\limits_{N}\sum\limits_{m=-\infty }^{\infty
}\int_{S^{1}}\int_{S^{1}}f(ank_{1}^{-1})\overline{%
f(a_{2}{}n_{2}k_{2})e^{-imk_{2}}}e^{-imk_{1}}dk_{1}dk_{2} \\
&&a^{-i\lambda }e^{-\text{ }i\langle \text{ }\xi ,\text{ }n\rangle
}a_{2}^{-i\lambda }e^{\text{ }i\langle \text{ }\xi ,\text{ }n_{2}\rangle
}dndadn_{2}da_{2}d\lambda d\xi \\
&=&\int\limits_{\mathbb{R}}\int\limits_{\mathbb{R}}\int\limits_{A}\int%
\limits_{N}\int\limits_{A}\int\limits_{N}\sum\limits_{m=-\infty }^{\infty
}\int_{S^{1}}\int_{S^{1}}f(ank_{1}^{-1})\overline{%
f(a_{2}{}n_{2}k_{2})e^{-imk_{2}}}e^{-imk_{1}}dk_{1}dk_{2} \\
&&a^{-i\lambda }e^{-\text{ }i\langle \text{ }\xi ,\text{ }n\rangle }%
\overline{a_{2}^{-i\lambda }e^{\text{ }-i\langle \text{ }\xi ,\text{ }%
n_{2}\rangle }}dndadn_{2}da_{2}d\lambda d\xi \\
&=&\int\limits_{\mathbb{R}}\int\limits_{\mathbb{R}}\sum\limits_{m=-\infty
}^{\infty }T\mathcal{F}f(\lambda ,\xi ,m)\overline{T\mathcal{F}f(\lambda
,\xi ,m)}d\lambda d\xi =\int\limits_{\mathbb{R}}\int\limits_{\mathbb{R}%
}\sum\limits_{m=-\infty }^{\infty }\left\vert T\mathcal{F}(f)(\lambda ,\xi
,m)\right\vert ^{2}d\lambda d\xi
\end{eqnarray*}

\section{Fourier Transform and Plancherel Formula $H.$}

\textbf{3.1.} Let $H$ be the real Heisenberg group of dimension $2n+1$ which
consists of all matrices of the form\newline
\begin{equation}
\left( 
\begin{array}{ccc}
1 & x & z \\ 
0 & I & y \\ 
0 & 0 & 1%
\end{array}%
\right) \text{ }
\end{equation}%
where $x\in \mathbb{R}^{n},$ $y\in \mathbb{R}^{n},$ $z\in \mathbb{R}$ and $%
\hspace{0.05in}I$ \hspace{0.05in}is the identity matrix of order $n$.

Let $H=\mathbb{R}^{n+1}\rtimes _{\iota }\mathbb{R}^{n}$ be the group of the
semi-direct product of the group $\mathbb{R}^{n+1}$ and $\mathbb{R}^{n}$,
via the group homomorphism\hspace{0.05in}$\iota :\mathbb{R}^{n}\rightarrow
Aut(\mathbb{R}^{n+1}),$ which is defined by:

\begin{equation}
\iota (x)(z,y)=(z+xy,y)=x(z,y)
\end{equation}%
for any $x=(x_{1},x_{2},$ $\ldots $ $,x_{n})\in \mathbb{R}^{n}$ , $%
y=(y_{1},y_{2},$ $\ldots $ $,y_{n})\in \mathbb{R}^{n}$ , $z\in \mathbb{R}$ ,
and $xy=\sum\limits_{i=1}^{n}x_{i}y_{i}$ , where $Aut(\mathbb{R}^{n+1})$ 
\hspace{0.05in}is the group of all automorphism of $\mathbb{R}^{n+1}$

\textbf{3.2.} Let $C^{\infty }(H)$ $,$ $\mathcal{D}(H)$ $,$ $\mathcal{D}%
^{\prime }(H)$ $,$ $\mathcal{E}^{\prime }(G)$ \hspace{0.05in}respectively
the space of $C^{\infty }$- functions , $C^{\infty }$ with compact support,
distribution and distribution with compact support on $G$. The Schwarts
space $\mathcal{S}(G)$ of $G$ can be considered as the Schwarts space $%
\mathcal{S}(\mathbb{R}^{2n+1})$ of the vector group $\mathbb{R}^{2n+1}$. The
action $\iota $ of the group $\mathbb{R}^{n}$ on $\mathbb{R}^{n+1}$ defines
a natural action $\iota $ on the dual $(\mathbb{R}^{n})^{\ast }$ of the
group $\mathbb{R}^{n+1}((\mathbb{R}^{n+1})^{\ast }\simeq \mathbb{R}^{n+1})$
which is given by :\newline
\begin{equation*}
x(\eta ,\lambda )=(\eta ,\eta x+\lambda )
\end{equation*}%
for any $\lambda \in \mathbb{R}^{n}$, $x\in \mathbb{R}^{n}$ and $\eta \in 
\mathbb{R}$ ,where ;\newline
\begin{equation*}
x(\eta ,\lambda )=\iota (x)(\eta ,\lambda )
\end{equation*}%
and 
\begin{equation*}
\eta x=\sum_{i=1}^{n}\eta x_{i}
\end{equation*}%
\textbf{Definition 3.1. }\textit{For every }$f\in $\textit{\ }$\mathcal{S}%
(G),$\textit{\ one can define its Fourier transform }$\mathcal{F}f$\textit{\ 
\hspace{0.05in}by :}

\begin{equation}
\mathcal{F}f(\xi )=\int\limits_{G}f(X)\text{ }e^{-i\left\langle \xi
,X\right\rangle }\text{ }dX
\end{equation}%
\textit{where} $X=((z,y);x)\in G,$ $\xi =((\eta ,\lambda );\mu )\in G,$ 
\textit{and} $dX=dz$ $dy$ $dx$ \textit{\hspace{0.05in}the Lebesgue measure on%
} $G$

\begin{equation*}
\left\langle \xi ,X\right\rangle =z\eta +y\lambda +x\mu =z\eta
+\sum\limits_{i=1}^n\lambda _iy_i+\sum\limits_{i=1}^nx_i\mu _i
\end{equation*}

\bigskip It is clear that the function $\mathcal{F}f\in \mathcal{S}(G)$ and
the mapping $f\mapsto \mathcal{F}f$ \hspace{0.05in}is a topological
isomorphism vector space $\mathcal{S}(G)$ onto it self.\bigskip

\medskip \textbf{Theorem 3.1.} \textit{The Fourier transform }$\mathcal{F}$%
\textit{\ satisfies :}

\begin{equation}
\overset{\vee }{g}\ast f(0)=\int\limits_{G}\mathcal{F}f(\xi )\text{ }%
\overline{\mathcal{F}g(\xi )\text{ }}d\xi
\end{equation}%
\textit{for every }$f\in \mathcal{S}(G)$\textit{\ and }$g\in \mathcal{S}(G)$%
\textit{, where }$\overset{\vee }{g}(X)=\overline{g(X^{-1})},$\textit{\ }$%
\xi =((\eta ,\lambda );\mu ),$\textit{\hspace{0.05in}}$d\xi =d\eta d\lambda
d\mu ,$\textit{\ is the Lebesgue measure on }$G=\mathbb{R}^{2n+1},$\textit{%
\hspace{0.05in}and }$\ast $\textit{\hspace{0.05in}denotes the convolution
product on }$G$\newline

\textit{Proof}\textbf{\ :} By the classical Fourier transform, we have:

\begin{eqnarray}
&&\overset{\vee }{g}\ast f(0)=\int\limits_{G}\mathcal{F}(\overset{\vee }{g}%
\ast f)(\xi )\text{ }d\xi  \notag \\
&=&\int\limits_{G}\int\limits_{G}\overset{\vee }{g}\ast f(X)\text{ }%
e^{-i\left\langle \xi ,X\right\rangle }\text{ }dX\text{ }d\xi  \notag \\
&=&\int\limits_{G}\int\limits_{G}\int\limits_{G}f(Y^{-1}X)\overline{g(Y^{-1})%
}\text{ }e^{-i\left\langle \xi ,X\right\rangle }\text{ }dY\text{ }dX\text{ }%
d\xi  \notag \\
&=&\int\limits_{G}\int\limits_{G}\int\limits_{G}f(YX)\overline{g(Y)}\text{ }%
e^{-i\left\langle \xi ,X\right\rangle }\text{ }dY\text{ }dX\text{ }d\xi
\end{eqnarray}

By change of variable $YX=X^{\prime },$with $X^{\prime }=((z,y);x)$ and $%
Y=((z^{\prime },y^{\prime });x^{\prime })$ we get :

\begin{eqnarray*}
X &=&Y^{-1}X^{\prime }=((-x^{\prime }(-z^{\prime },-y^{\prime }))-x^{\prime
})((z,y);x) \\
\text{ \hspace{0.05in}} &=&((-x^{\prime }(z-z^{\prime },y-y^{\prime
}));x-x^{\prime })
\end{eqnarray*}%
this gives us :

\begin{eqnarray}
&&e^{-i\left\langle \xi ,X\right\rangle }  \notag \\
&=&e^{-i\left\langle \xi ,Y^{-1}X^{\prime }\right\rangle }  \notag \\
&=&e^{-i\left\langle (-x^{\prime }(\eta ,\lambda );\mu );((z-z^{\prime
},y-y^{\prime });x-x^{\prime })\right\rangle }  \notag \\
&=&e^{-i\left\langle ((\eta ,-\eta x^{\prime }+\lambda );\mu ),((z-z^{\prime
},y-y^{\prime });x-x^{\prime })\right\rangle }
\end{eqnarray}

By the invariant of the Lebesgue measures $d\eta ,$ $d\lambda ,$ and $d\mu $
we obtain,

\begin{eqnarray}
&&\overset{\vee }{g}\ast f(0)=\int\limits_{G}f(X)\text{ }\overline{g(X)\text{
}}dX  \notag \\
&=&\int\limits_{G}\int\limits_{G}\int\limits_{G}f(X)e^{-i\left\langle \xi
,X\right\rangle }\text{ }\overline{g(Y)e^{-i\left\langle \xi ,Y\right\rangle
}}\text{ }\hspace{0.05in}dX\text{ }dY\hspace{0.05in}d\xi  \notag \\
&=&\int\limits_{G}\mathcal{F}f(\xi )\text{ }\overline{\mathcal{F}g(\xi )}%
\text{ }d\xi
\end{eqnarray}%
where $0=((0,0);0)$ \hspace{0.05in}is the identity of $G$ , whence the
theorem.

\textbf{Corollary 3.1. }\textit{In theorem \textbf{3.1} , if\hspace{0.05in}%
we take }$\hspace{0.05in}g=\overset{\vee }{f},$\textit{\ we obtain the
Plancherel formula on }$G$

\begin{equation}
\overset{\vee }{f}\ast f(0)=\int\limits_{G}\left\vert f(X)\right\vert
^{2}dX=\int\limits_{\mathbb{R}^{n}}\left\vert \mathcal{F}f(\xi )\right\vert
^{2}d\xi
\end{equation}

\section{Fourier Transform and Plancherel Formula on the Real Jacobi group $%
N\rtimes SL(2,\mathbb{R})$}

\bigskip \textbf{4.1. }Let $H$ be the $3-$dimensional Heisenberg group, with
multiplication%
\begin{equation}
(z_{1},y_{1},x_{1})(z_{2},y_{2},x_{2})=(z_{1}+z_{2}+x_{1}y_{2}-x_{2}y_{1},y_{1}+x_{1},x_{2}+y_{2})
\end{equation}

\bigskip The group $H$ is isomorphic onto the following Heisenberg group of
all matrices 
\begin{equation}
N=\{\mathbb{R}^{2}\rtimes _{\sigma }\mathbb{R}\simeq \left( 
\begin{array}{ccc}
1 & x & z \\ 
0 & 1 & z \\ 
0 & 0 & 1%
\end{array}%
\right) :\text{ }(z,y,x)\in \mathbb{R}^{3}\text{ }\}
\end{equation}%
where $\sigma :\mathbb{R\rightarrow }Aut\mathbb{\ (R}^{2})$ is the group
homomorphism from the real group into the group $Aut\mathbb{\ (R}^{2})$ of
all automorphisms of the vector group $\mathbb{R}^{2},$ defined as 
\begin{equation*}
\sigma (x)(z,y)=(z+xy,y)
\end{equation*}

\bigskip So the group $H$ can be identified with the group $N$, where the
multiplication becomes as%
\begin{equation}
(z_{1},y_{1},x_{1})(z_{2},y_{2},x_{2})=(z_{1}+z_{2}+x_{1}y_{2},y_{1}+x_{1},x_{2}+y_{2})
\end{equation}

\bigskip Now we define the Jacobi group $J$ as $N\rtimes _{\rho }SL(2,%
\mathbb{R})$ the semidirect of the Heisenberg group $N$ and the real
semisimple Lie group $SL(2,\mathbb{R}),$ where $\rho :SL(2,\mathbb{R}%
)\rightarrow Aut\mathbb{\ (}N)$ is the group homomorphism from the real
group into the group $Aut\mathbb{\ (}N)$ of all automorphisms of the vector
group $N,$ defined as 
\begin{eqnarray*}
\rho (M)(z,y,x) &=&(z,\left[ 
\begin{array}{cc}
y & x%
\end{array}%
\right] M) \\
&=&(z,\left[ 
\begin{array}{cc}
y & x%
\end{array}%
\right] \left[ 
\begin{array}{cc}
a & b \\ 
c & d%
\end{array}%
\right] ) \\
&=&(z,ya+xc,yb+xd)
\end{eqnarray*}%
where $M=\left[ 
\begin{array}{cc}
a & b \\ 
c & d%
\end{array}%
\right] $. So, any elemet $g\in J$ can be written in an unique way as $g=(X,$
$M)$ with $M\in SL(2,\mathbb{R})$ and $X=(z,y,x)\in N.$ Multiplication in $J$
is then given as 
\begin{eqnarray}
&&(X_{1},M_{1})(X_{2},M_{2})  \notag \\
&=&(z_{1},y_{1},x_{1},M_{1})(z_{2},y_{2},x_{2},M_{2}) \\
&=&((z_{1},y_{1},x_{1})(z_{2,}\left[ 
\begin{array}{cc}
y_{2} & x_{2}%
\end{array}%
\right] \left[ 
\begin{array}{cc}
a_{1} & b_{1} \\ 
c_{1} & d_{1}%
\end{array}%
\right] ),\left[ 
\begin{array}{cc}
a_{1} & b_{1} \\ 
c_{1} & d_{1}%
\end{array}%
\right] \left[ 
\begin{array}{cc}
a_{2} & b_{2} \\ 
c_{2} & d_{2}%
\end{array}%
\right] )  \notag \\
&=&((z_{1},y_{1},x_{1})(z_{2,}y_{2}a_{1}+x_{2}c_{1},y_{2}b_{1}+x_{2}d_{1}),M_{1}M2)
\notag \\
&=&(z_{1}+z_{2}+x_{1}y_{2}a_{1}+x_{1}x_{2}c_{1},y_{1}+y_{2}a_{1}+x_{2}c_{1},x_{1}+y_{2}b_{1}+x_{2}d_{1},M_{1}M2)
\end{eqnarray}%
where $X_{1}=(z_{1},y_{1},x_{1})\in N,X_{2}=(z_{2},y_{2},x_{2})\in N,M_{1}=$ 
$\left[ 
\begin{array}{cc}
a_{1} & b_{1} \\ 
c_{1} & d_{1}%
\end{array}%
\right] \in SL(2,\mathbb{R}),$and $M_{2}=\left[ 
\begin{array}{cc}
a_{1} & b_{1} \\ 
c_{1} & d_{1}%
\end{array}%
\right] \in SL(2,\mathbb{R})$

From now on, our useful for the multiplication in $J$ will be as 
\begin{equation*}
(X_{1},M_{1})(X_{2},M_{2})=(X_{1}\rho (M_{1})(X_{2}),M_{1}M_{2})
\end{equation*}

\textbf{Definition 4.1}. \textit{Let }$Q=H\times SL(2,\mathbb{R})\rtimes
_{\rho }SL(2,\mathbb{R})$ \textit{be the group with law}: 
\begin{eqnarray}
X\cdot Y &=&(X_{1},M_{1},M_{2})(Y_{1},N_{1},N_{2})  \notag \\
&=&(X_{1}\rho (M_{2})(Y_{1}),M_{1}+N_{1},M_{2}+N_{2})
\end{eqnarray}%
\textit{for all} $X=(X_{1},M_{1},M_{2})$ $\in Q$ \textit{and }$%
Y=(Y_{1},N_{1},N_{2})\in Q.$ From defintion \textbf{4.1. }the Jacobi group $%
J $ can be identified with a subgroup $N\times \{I_{SL(2,\mathbb{R}%
}\}\rtimes _{\rho }SL(2,\mathbb{R})$ of $Q\ .$ Let $A=N\times SL(2,\mathbb{R}%
)\times _{\rho }\{I_{SL(2,\mathbb{R}}\}$ be the subgroup of $Q,$ which is
the direct product of $N$ with $SL(2,\mathbb{R})$

\textbf{Definition 4.2}. \textit{For any function} $f\in \mathcal{D}(J),$ 
\textit{we can define a function} $\widetilde{f}$ \textit{on }$Q$ \textit{by}%
\begin{equation}
\widetilde{f}(X,M_{1},M_{2})=f(MX,M_{1}M_{2})
\end{equation}

\textbf{Remark 4.1. }The function $\widetilde{f}$ \ is invariant in the
following sense

\begin{equation}
\widetilde{f}(N^{-1}X,M_{1},N^{-1}M_{2})=\widetilde{f}%
(N^{-1}X,M_{1},N^{-1}M_{2})
\end{equation}

\textbf{Theorem 4.1. }\textit{For any function }$\psi \in \mathcal{D}(J)$ 
\textit{and }$\widetilde{f}\in \mathcal{D}(Q)$ \textit{invariant in sense }$%
(32)$\textit{, we get}%
\begin{equation}
\psi \ast \widetilde{f}(X,M_{1},M_{2})=\widetilde{f}\ast _{c}\psi
(X,M_{1},M_{2})
\end{equation}%
\textit{\ where} $\ast $ \textit{signifies the convolution product on} $J$ 
\textit{with respect the variable} $(X,M_{2}),$ \textit{and }$\ast _{c}$%
\textit{signifies the convolution product on} $A$ \textit{with respect the
variable} $(X,M_{1})$

\textit{Proof : }In fact we have%
\begin{eqnarray}
&&\psi \ast \widetilde{f}(X,M_{1},M_{2})  \notag \\
&=&\dint\limits_{N}\int_{SL(2,\mathbb{R)}}\widetilde{f}%
((Y,M)^{-1}(X,M_{1},M_{2}))\psi (Y,M)dYdM  \notag \\
&=&\dint\limits_{N}\int_{SL(2,\mathbb{R)}}\widetilde{f}%
[(M{}^{-1}(-Y),M{}^{-1})(X,M_{1},M_{2})]\psi (Y,M)dYdM  \notag \\
&=&\dint\limits_{N}\int_{SL(2,\mathbb{R)}}\widetilde{f}%
[(M{}^{-1}(Y-X),M_{1},M^{-1}M_{2})(X,M_{1},M_{2})]\psi (Y,M)dYdM  \notag \\
&=&\dint\limits_{N}\int_{SL(2,\mathbb{R)}}\widetilde{f}%
[(M^{-1}(Y-X),M_{1},M^{-1}M_{2})]\psi (v^{\prime },g^{\prime })dYdM  \notag
\\
&=&\dint\limits_{N}\int_{SL(2,\mathbb{R)}}\widetilde{f}%
[Y-X,M_{1}M{}^{-1},M_{2}]\psi (Y,M)dYdM=\widetilde{f}\ast _{c}\psi
(X,M_{1},M_{2})
\end{eqnarray}%
\textbf{\ }for any function $\psi \in \mathcal{D}(J)$ and\textit{\ }$%
\widetilde{f}\in \mathcal{D}(Q)$

\bigskip \textbf{Definition 4}.\textbf{3}. \textit{For any }$k_{1}$ $\in
S^{1}$ \textit{let }$\Gamma _{k_{1}}\Psi $ \textit{be the fuction defined by}%
\begin{equation}
\Gamma _{k_{1}}\Psi (v,g)=\Psi (X,gk_{1})
\end{equation}%
\textit{for any} $v\in N$, \ $g\in SL(2,\mathbb{R)}$ \textit{and }$k_{1}$ $%
\in S^{1}$

\textbf{Definition 4.4.}\textit{\ Let} $f\in C_{0}^{\infty }(J),$ \textit{we
define its Fourier transform by}%
\begin{equation}
\mathcal{F}_{N}T\mathcal{F}\Psi (\eta ,m,\xi ,\lambda
)=\dint\limits_{N}\int_{A}\int_{N}\int_{S^{1}}\Psi (v,kna)e^{-\text{ }%
i\langle \text{ }\eta ,\text{ }v\rangle }\text{ }e^{-ikm}a^{-i\lambda }e^{-%
\text{ }i\langle \text{ }\xi ,\text{ }n\rangle }dkdadndv  \notag
\end{equation}%
\textit{where} $\mathcal{F}_{N}$ \textit{is the Fourier transform on} $N,$ $%
kna=g,$ $\eta =(\eta _{1},\eta _{2},\eta _{3})\in \mathbb{R}^{3},$ $%
v=(v_{1},v_{2},v_{3})\in N,$ \textit{and} $dv=dv_{1}dv_{2}dv_{3}$ \textit{is
the Lebesgue measure on} $N$

\begin{eqnarray}
\langle \eta ,v\rangle &=&\langle (\eta _{1},\eta _{2},\eta
_{3}),(v_{1},v_{2},v_{3})\rangle  \notag \\
&=&\eta _{1}v_{1}+\eta _{2}v_{2}+v_{3}\eta _{3}
\end{eqnarray}

\textbf{\textit{Plancherel's Theorem }4.2\textit{. }}\textit{For any
function }$f\in $\textit{\ }$L^{1}(J)\cap $\textit{\ }$L^{2}(J),$\textit{we
get}%
\begin{equation}
\int_{J}\left\vert \Psi (v,g)\right\vert
^{2}dvdg=\dint\limits_{N}\int\limits_{\mathbb{R}}\int\limits_{\mathbb{R}%
}\sum\limits_{m=-\infty }^{\infty }\left\vert \mathcal{F}_{H}T\mathcal{F}%
\Psi (\eta ,m,\xi ,\lambda )\right\vert d\eta d\lambda d\xi
\end{equation}

\bigskip Proof: For any function $\Psi \in \mathit{\ }L^{1}(J)\cap \mathit{\ 
}L^{2}(J),$ we get \ \ 
\begin{eqnarray*}
&&\Gamma _{I_{K}}\Psi \ast \widetilde{\overset{\vee }{\Psi }}(0,I_{SL(2,%
\mathbb{R)}},I_{SL(2,\mathbb{R)}}) \\
&=&\sum\limits_{m=-\infty }^{\infty }[\Gamma _{k_{1}}\Psi \ast \widetilde{%
\overset{\vee }{\Psi }}(0,I_{SL(2,\mathbb{R)}},I_{SL(2,\mathbb{R)}%
})e^{-ikm}dk_{1}] \\
&=&\dint\limits_{N}\dint_{SL(2,\mathbb{R})}\sum\limits_{m=-\infty }^{\infty
}[\dint_{S^{1}}\widetilde{\overset{\vee }{\Psi }}((w,g{})^{-1}(0,I_{SL(2,%
\mathbb{R)}},I_{SL(2,\mathbb{R)}}))\Gamma _{k_{1}}\Psi
(w,g)e^{-ikm}dk_{1}]dwdg \\
&=&\dint\limits_{N}\dint_{SL(2,\mathbb{R})}\sum\limits_{m=-\infty }^{\infty
}[\dint_{S^{1}}\widetilde{\overset{\vee }{\Psi }}%
(g{}^{-1}(0-w),I_{G},g^{-1}I_{G})\Psi (w,k_{1}g{})e^{-ikm}dk_{1}]dwdg \\
&=&\dint\limits_{N}\dint_{SL(2,\mathbb{R})}\sum\limits_{m=-\infty }^{\infty
}[\dint_{S^{1}}\widetilde{\overset{\vee }{\Psi }}(-w,I_{G}g^{-1},I_{G})\Psi
(w,k_{1}g{})e^{-ikm}dk_{1}]dwdg \\
&=&\dint\limits_{N}\dint_{SL(2,\mathbb{R})}\sum\limits_{m=-\infty }^{\infty
}[\dint_{S^{1}}\overset{\vee }{\Psi }(g^{-1}(-w),g^{-1})\Psi
(w,k_{1}g{})e^{-ikm}dk_{1}]dwdg \\
&=&\dint\limits_{N}\dint_{SL(2,\mathbb{R})}\overset{\vee }{\Psi }%
(g^{-1}(-w),g^{-1})\Psi (w,I_{S^{1}}g{})]dwdg \\
&=&\dint\limits_{H}\dint_{SL(2,\mathbb{R})}\overline{\Psi \lbrack
(g^{-1}(-w),g^{-1})^{-1}]}\Psi (w,g{})]dwdg \\
&=&\dint\limits_{H}\dint_{SL(2,\mathbb{R})}\overline{\Psi (w,g)}\Psi
(w,g{})]dwdg \\
&=&\dint\limits_{H}\dint_{SL(2,\mathbb{R})}\left\vert \Psi (w,g)\right\vert
^{2}dwdg=\dint_{J}\left\vert \Psi (w,g)\right\vert ^{2}dwdg
\end{eqnarray*}

In other hand

\begin{eqnarray*}
&&\Gamma _{I_{K}}\Psi \ast \widetilde{\overset{\vee }{\Psi }}(0,I_{G},I_{G})
\\
&=&\dint\limits_{H}\dint\limits_{\mathbb{R}}\int\limits_{\mathbb{R}%
}\sum\limits_{m=-\infty }^{\infty }[\mathcal{F}_{H}T\mathcal{F(}\Gamma
_{k_{1}}\Psi \ast \widetilde{\overset{\vee }{\Psi }})(\eta ,\xi ,\lambda
,m,I_{G})e^{-ik_{1}m}dk_{1}]d\eta d\xi d\lambda \\
&=&\dint\limits_{H}\dint\limits_{\mathbb{R}}\int\limits_{\mathbb{R}%
}\sum\limits_{m=-\infty }^{\infty }[\mathcal{F}_{H}T\mathcal{F}[\Gamma
_{k_{1}}\Psi \ast \widetilde{\overset{\vee }{\Psi }}(v,I_{k}na,I_{G})\gamma
(k_{1}^{-1})dk_{1}] \\
&&e^{-i\langle v,\eta \rangle }e^{-i\langle n,\xi \rangle }a^{-i\lambda
}dvdndad\eta d\xi d\lambda \\
&=&\dint\limits_{H}\dint\limits_{\mathbb{R}}\int\limits_{\mathbb{R}%
}\sum\limits_{m=-\infty }^{\infty }[\mathcal{F}_{H}T\mathcal{F}%
[\dint\limits_{K}\widetilde{\overset{\vee }{\Psi }}%
((w,g_{2}{})^{-1}(v,I_{k}na,I_{G})\Gamma _{k_{1}}\Psi (w,g_{2}) \\
&&\gamma (k_{1}^{-1}))dk_{1}]e^{-i\langle v,\eta \rangle }e^{-i\langle n,\xi
\rangle }a^{-i\lambda }dvdndadwdg_{2}d\eta d\xi d\lambda \\
&=&\dint\limits_{H}\dint\limits_{\mathbb{R}}\int\limits_{\mathbb{R}%
}\sum\limits_{m=-\infty }^{\infty }[\mathcal{F}_{H}T\mathcal{F}%
[\dint\limits_{K}\widetilde{\overset{\vee }{\Psi }}%
((g_{2}{}^{-1}(-w),g_{2}{}^{-1})(v,I_{k}na,I_{G})\Gamma _{k_{1}}\Psi
(w,g_{2}) \\
&&\gamma (k_{1}^{-1}))dk_{1}]e^{-i\langle v,\eta \rangle }e^{-i\langle n,\xi
\rangle }a^{-i\lambda }dvdndadwdg_{2}d\eta d\xi d\lambda \\
&=&\dint\limits_{G}\dint\limits_{\mathbb{R}^{4}}\dint\limits_{\mathbb{R}%
}\int\limits_{\mathbb{R}^{2}}\dint\limits_{\mathbb{R}^{4}}\dint\limits_{%
\mathbb{R}}\int\limits_{\mathbb{R}^{2}}\dint\limits_{\mathbb{R}%
^{4}}\dsum\limits_{\gamma \in \widehat{K}}d_{\gamma }tr[\dint\limits_{K}%
\widetilde{\overset{\vee }{\Psi }}%
((g_{2}{}^{-1}(-w)+(g_{2}{}^{-1})(v),I_{k}na,g_{2}{}^{-1}I_{G}) \\
&&\Gamma _{k_{1}}\Psi (w,g_{2})\gamma (k_{1}^{-1}))dk_{1}]e^{-i\langle
v,\eta \rangle }e^{-i\langle n,\xi \rangle }a^{-i\lambda
}dvdndadwdg_{2}d\eta d\xi d\lambda \\
&=&\dint\limits_{N}\dint\limits_{\mathbb{R}}\int\limits_{\mathbb{R}%
}\sum\limits_{m=-\infty }^{\infty }[\mathcal{F}_{H}T\mathcal{F}%
[\dint\limits_{K}\widetilde{\overset{\vee }{\Psi }}%
((g_{2}{}^{-1}(v-w),I_{k}na,g_{2}{}^{-1}I_{G}) \\
&&\Psi (w,k_{1}g_{2})\gamma (k_{1}^{-1}))dk_{1}]e^{-i\langle v,\eta \rangle
}e^{-i\langle n,\xi \rangle }a^{-i\lambda }dvdndadwdg_{2}d\eta d\xi d\lambda
\\
&=&\dint\limits_{N}\dint\limits_{\mathbb{R}}\int\limits_{\mathbb{R}%
}\sum\limits_{m=-\infty }^{\infty }[\mathcal{F}_{H}T\mathcal{F}%
[\dint\limits_{K}\widetilde{\overset{\vee }{\Psi }}%
(((v-w),nag_{2}^{-1}{},I_{G})\Psi (w,k_{1}g_{2}) \\
&&\gamma (k_{1}^{-1}))dk_{1}]e^{-i\langle v,\eta \rangle }e^{-i\langle n,\xi
\rangle }a^{-i\lambda }dvdndadwdg_{2}d\eta d\xi d\lambda \\
&=&\dint\limits_{N}\dint\limits_{\mathbb{R}}\int\limits_{\mathbb{R}%
}\sum\limits_{m=-\infty }^{\infty }[\mathcal{F}_{H}T\mathcal{F}%
[\dint\limits_{K}\dint\limits_{K}\widetilde{\overset{\vee }{\Psi }}%
((v-w,I_{k}naa_{2}{}^{-1}n_{2}{}^{-1}k_{2}{}^{-1},I_{G})\Psi
(w,k_{1}k_{2}n_{2}a_{2}) \\
&&\gamma (k_{1}^{-1}))dk_{1}dk_{2}]e^{-i\langle v,\eta \rangle }e^{-i\langle
n,\xi \rangle }a^{-i\lambda }dvdndadwdn_{2}da_{2}d\eta d\xi d\lambda
\end{eqnarray*}

We continue our calculation$.$%
\begin{eqnarray*}
&&\Gamma _{I_{K}}\Psi \ast \widetilde{\overset{\vee }{\Psi }}(0,I_{G},I_{G})
\\
&=&\dint\limits_{N}\dint\limits_{\mathbb{R}}\int\limits_{\mathbb{R}%
}\sum\limits_{m=-\infty }^{\infty }[\mathcal{F}_{H}T\mathcal{F}%
[\dint\limits_{K}\dint\limits_{K}\widetilde{\overset{\vee }{\Psi }}%
((v,ank_{2}{}^{-1},I_{G})\Psi (w,k_{1}k_{2}n_{2}a_{2}) \\
&&\gamma (k_{1}^{-1}))dk_{1}dk_{2}]e^{-i\langle v,\eta \rangle }e^{-i\langle
n,\xi \rangle }a^{-i\lambda }e^{-i\langle w,\eta \rangle }e^{-i\langle
n_{2},\xi \rangle }a_{2}^{-i\lambda }dvdndadwdn_{2}da_{2}d\eta d\xi d\lambda
\\
&=&\dint\limits_{N}\dint\limits_{\mathbb{R}}\int\limits_{\mathbb{R}%
}\sum\limits_{m=-\infty }^{\infty }[\mathcal{F}_{H}T\mathcal{F}%
[\dint\limits_{K}\dint\limits_{K}\widetilde{\overset{\vee }{\Psi }}%
((v,ank_{2}{}^{-1}k_{1},I_{G})\Psi (w,k_{2}n_{2}a_{2})\gamma
(k_{1}^{-1})\gamma (k_{2}^{-1})) \\
&&dk_{1}dk_{2}]e^{-i\langle v,\eta \rangle }e^{-i\langle n,\xi \rangle
}a^{-i\lambda }e^{-i\langle w,\eta \rangle }e^{-i\langle n_{2},\xi \rangle
}a_{2}^{-i\lambda }dvdndadwdn_{2}da_{2}d\eta d\xi d\lambda \\
&=&\dint\limits_{N}\dint\limits_{\mathbb{R}}\int\limits_{\mathbb{R}%
}\sum\limits_{m=-\infty }^{\infty }[\mathcal{F}_{H}T\mathcal{F}%
[\dint\limits_{K}\dint\limits_{K}\widetilde{\overset{\vee }{\Psi }}%
((v,ank_{1},I_{G})\Psi (w,k_{2}n_{2}a_{2})\gamma (k_{1}^{-1})\gamma
(k_{2}^{-1})) \\
&&dk_{1}dk_{2}]e^{-i\langle v,\eta \rangle }e^{-i\langle n,\xi \rangle
}a^{-i\lambda }e^{-i\langle w,\eta \rangle }e^{-i\langle n_{2},\xi \rangle
}a_{2}^{-i\lambda }dvdndadwdn_{2}da_{2}d\eta d\xi d\lambda \\
&=&\dint\limits_{N}\dint\limits_{\mathbb{R}}\int\limits_{\mathbb{R}%
}\sum\limits_{m=-\infty }^{\infty }[\mathcal{F}_{H}T\mathcal{F}%
[\dint\limits_{K}\dint\limits_{K}\overset{\vee }{\Psi }(ank_{1}v,ank_{1})%
\Psi (w,k_{2}n_{2}a_{2})\gamma (k_{1}^{-1})\gamma (k_{2}^{-1})) \\
&&dk_{1}dk_{2}]e^{-i\langle v,\eta \rangle }e^{-i\langle n,\xi \rangle
}a^{-i\lambda }e^{-i\langle w,\eta \rangle }e^{-i\langle n_{2},\xi \rangle
}a_{2}^{-i\lambda }dvdndadwdn_{2}da_{2}d\eta d\xi d\lambda \\
&=&\dint\limits_{H}\dint\limits_{\mathbb{R}}\int\limits_{\mathbb{R}%
}\sum\limits_{m=-\infty }^{\infty }[\mathcal{F}_{H}T\mathcal{F}%
[\dint\limits_{K}\dint\limits_{K}\overline{\Psi ((ank_{1}(v),ank_{1})^{-1})}%
\Psi (w,k_{2}n_{2}a_{2})\gamma (k_{1}^{-1}) \\
&&\gamma (k_{2}^{-1}))dk_{1}]e^{-i\langle v,\eta \rangle }e^{-i\langle n,\xi
\rangle }a^{-i\lambda }e^{-i\langle w,\eta \rangle }e^{-i\langle n_{2},\xi
\rangle }a_{2}^{-i\lambda }dvdndadwdk_{2}dn_{2}da_{2}d\eta d\xi d\lambda \\
&=&\dint\limits_{H}\dint\limits_{\mathbb{R}}\int\limits_{\mathbb{R}%
}\sum\limits_{m=-\infty }^{\infty }[\mathcal{F}_{H}T\mathcal{F}%
[\dint\limits_{K}\dint\limits_{K}\overline{\Psi (-v,k_{1}^{-1}n^{-1}a^{-1})}%
\Psi (w,k_{2}n_{2}a_{2})\gamma (k_{1}^{-1})\gamma (k_{2}^{-1})) \\
&&dk_{1}dk_{2}]e^{-i\langle v,\eta \rangle }e^{-i\langle n,\xi \rangle
}a^{-i\lambda }e^{-i\langle w,\eta \rangle }e^{-i\langle n_{2},\xi \rangle
}a_{2}^{-i\lambda }dvdndadwdn_{2}da_{2}d\eta d\xi d\lambda \\
&=&\dint\limits_{H}\dint\limits_{\mathbb{R}}\int\limits_{\mathbb{R}%
}\sum\limits_{m=-\infty }^{\infty }[\mathcal{F}_{H}T\mathcal{F}%
[\dint\limits_{K}\dint\limits_{K}\overline{\Psi (v,k_{1}na)}\Psi
(w,k_{2}n_{2}a_{2})\gamma ^{\ast }(k_{1}^{-1})\gamma (k_{2}^{-1})) \\
&&dk_{1}dk_{2}]\overline{e^{-i\langle v,\eta \rangle }e^{-i\langle n,\xi
\rangle }a^{-i\lambda }}e^{-i\langle w,\eta \rangle }e^{-i\langle n_{2},\xi
\rangle }a_{2}^{-i\lambda }dvdndadwdn_{2}da_{2}d\eta d\xi d\lambda \\
&=&\dint\limits_{\mathbb{R}}\int\limits_{\mathbb{R}^{2}}\dint\limits_{%
\mathbb{R}^{4}}\dsum\limits_{\gamma \in \widehat{K}}d_{\gamma }tr[\mathcal{F}%
_{\mathbb{R}^{4}}T\mathcal{F}\overline{\Psi (\eta ,\gamma ^{\ast },\xi
,\lambda {})}\mathcal{F}_{\mathbb{R}^{4}}T\mathcal{F}\Psi (\eta ,\gamma ,\xi
,\lambda )]d\eta d\xi d\lambda \\
&=&\dint\limits_{\mathbb{R}}\int\limits_{\mathbb{R}^{2}}\dint\limits_{%
\mathbb{R}^{4}}\dsum\limits_{\gamma \in \widehat{K}}d_{\gamma }\left\Vert 
\mathcal{F}_{\mathbb{R}^{4}}T\mathcal{F}\Psi (\eta ,\gamma ,\xi ,\lambda
)\right\Vert _{H.S}^{2}d\eta d\xi d\lambda
\end{eqnarray*}

\section{\protect\bigskip Left Ideals of the Group Algebra $L^{1}(N).$}

First, I will prove the solvability of any invariant differential operator
on the connected solvable group $N=\mathbb{R}^{2}\rtimes _{\sigma }\mathbb{R}%
.$ Therefor, I will extend the group by a larger group $E=\mathbb{R}%
^{2}\times \mathbb{R}\times \mathbb{R}$, with multiplication $(n,a,b)$ and $%
(m,x,y)$ as

\begin{equation}
(n,a,b)(m,x,y)=(n+m+\sigma (b)x,a+x,b+y)
\end{equation}

\bigskip Let $\textfranc =\mathbb{R}^{2}\times \mathbb{R}$ be the abelian
group, which is the direct product of two real vector groups $\mathbb{R}^{2}$
and $\mathbb{R}$

\textbf{Definition 5.1}. \textit{For any function} $f\in \mathcal{D}(N),$ 
\textit{we can define a function} $\tau f$ \textit{on }$E$ \textit{by}%
\begin{equation}
\tau f(n,a,b)=f(\sigma (a)n,ab)
\end{equation}

\textbf{Remark 5.1. }The function $\tau f$ \ is invariant in the following
sense

\begin{equation}
\tau f(\sigma (x^{-1})n,xa,x^{-1}b)=\tau f(n,a,b)
\end{equation}%
Therefor denote by $\tau C^{\infty }(N)$ $(resp.$ $\tau C^{\infty
}(\textfranc )$ ) the image of $C^{\infty }(N$ $)$ $(resp.C^{\infty
}(\textfranc )$\ then we have%
\begin{equation*}
\ \tau C^{\infty }(N)|_{N}=C^{\infty }(N)
\end{equation*}%
\begin{equation}
\ \tau C^{\infty }(\textfranc )|_{\textfranc }=C^{\infty }(\textfranc )
\end{equation}

\bigskip \textbf{Definition 5.2. }\textit{Let be the mapping} $\Lambda :\tau
C^{\infty }(E)|_{\textfranc \text{ }}$ $\longrightarrow \tau C^{\infty
}(E)|_{N\text{ }}$\textit{\ defined by}%
\begin{equation}
\Lambda (\tau f|_{\textfranc }\ )(z,y,0)=\tau f|_{N}(z,0,y)
\end{equation}%
\textit{is topological isomorphisms and its inverse is nothing but }$\Gamma
^{-1}$\textit{defined by}

\begin{equation}
\Lambda ^{-1}(\tau f|_{N}\ )(n,0,a)=\tau f|_{\textfranc }(n,a,0)
\end{equation}

My main result is

\textbf{Theorem 5.1.} \textit{If }$P_{u}$\ \textit{any} \textit{invariant
differential operator on }$N$\textit{\ associated to the distribution }$u\in 
\mathcal{U}$\textit{, then, we have }

\begin{equation}
P_{u}\text{ }C^{\infty }(N)=C^{\infty }(N)
\end{equation}

\textit{Proof:} Let $Q_{u}$\ be the invariant differential operator with
constant coefficients on $K$ associated to $u$ , then by the theory of
differential operators with constant coefficients $[20]$, we get 
\begin{equation}
Q_{u}\ \ \tau C^{\infty }(E)|_{\textfranc \text{ }}=\tau C^{\infty
}(E)|_{\textfranc \text{ }}=C^{\infty }(\textfranc )
\end{equation}

That means for any $\psi (n,a)\in C^{\infty }(\textfranc ),$ there exist a
function $\varphi (n,a,x)\in \tau C^{\infty }(E)|_{\textfranc \text{ }},$
such that%
\begin{equation}
Q_{u}\varphi (n,a,1)=u\ast _{c}\varphi (n,a,0)=\psi (n,a)
\end{equation}

The function $\psi (n,a)$ can be \ transformed as an invariant function $%
\psi \in \tau C^{\infty }(E)|_{\textfranc \text{ }}$ as follows

\begin{equation}
\psi (n,a)=\tau \psi (\rho (a^{-1})n,a,0)
\end{equation}

In other side, we have%
\begin{eqnarray}
&&\Lambda Q_{u}\ \ \varphi (n,a,0)  \notag \\
&=&Q_{u}\ \ \varphi (n,0,a)=u\ast _{c}\varphi (n,0,a)  \notag \\
&=&u\ast \varphi (n,1,a)=P_{u}\text{ }\varphi (n,0,a)  \notag \\
&=&\Lambda \tau \psi (\rho (a^{-1})n,a,1)=\tau \psi (\sigma (a^{-1})n,1,a) 
\notag \\
&=&\psi (n,a)
\end{eqnarray}

So the proof of the solvability of any right invariant differential operator
on $N.$

If $I$ is a subspace of $L^{1}(N),$ we denote $\tau I$ its image by the
mapping $\tau $, let $J=\tau $ $I|_{\textfranc }.$ My main result is:

\textbf{Theorem 5.2.} \textit{Let} $I$ \textit{be a subspace of} $L^{1}(N),$ 
\textit{then the following conditions are equivalents}.

$(i)$ $J=\tau I|_{\textfranc }$ \textit{is an ideal in the Banach algebra} $%
L^{1}(\textfranc ).$

$(ii)$ $I$ \textit{is a left ideal in the Banach algebra} $L^{1}(N).$

\bigskip \textit{Proof:} $(i)$ implies $(ii)$\ Let $I$ be a subspace of the
space $L^{1}(N)$ and $\tau I$ the image of $I$ by $\tau $ such that $J=\tau
I|_{\textfranc }$ \ is an ideal in $L^{1}(\textfranc ),$ then we have: 
\begin{equation}
u\ast _{c}\tau I|_{\textfranc }(n,a,0)\subseteq \tau I|_{\textfranc }(n,a,0)
\end{equation}%
for any $u\in L^{1}(\textfranc )$ and $(n,a)\in \textfranc $, where%
\begin{equation}
u\ast _{c}\tau I|_{\textfranc }(n,a,0)=\left\{ \int\limits_{\textfranc }\tau
f|_{\textfranc }\ \left[ n-m,a-b,0)\right] u(m,b)dm\frac{db}{b},\text{ }f\in 
\text{ }I\right\}
\end{equation}

It shows that%
\begin{equation}
u\ast _{c}\tau f|_{\textfranc }\ (n,a,0)\in \tau I|_{\textfranc }(n,a,0)
\end{equation}%
for any $\tau f\in \tau I.$ Apply equation$(32),$ we get%
\begin{eqnarray}
&&\Gamma (u\ast _{c}\tau f|_{\textfranc }\ )(n,a,0)  \notag \\
&=&u\text{ }\ast \tau f(n,1,a)\in \Gamma (\tau I|_{\textfranc }(n,a,0) 
\notag \\
&=&\tau I|_{N}(n,0,a)=I
\end{eqnarray}%
$(ii)$ implies $(i),$ if $I$ is an ideal in $L^{1}(N),$ then we get 
\begin{eqnarray}
&&u\ast \tau I|_{N}\ (n,0,a)  \notag \\
&=&u\ast I\ (n,a)\subseteq \tau I\ |_{N}(n,0,a)=I\ (n,a)
\end{eqnarray}%
where%
\begin{equation}
u\ast \tau I\ |_{N}(n,1,a)=\left\{ \int\limits_{N}\tau f|_{N}\ \left[ \sigma
(-b)(n-m),1,a-b\right] u(m,b)dm\frac{db}{b},\text{ }f\in \text{ }I\right\}
\end{equation}

By equation $(36),$we obtain%
\begin{eqnarray}
&&\chi ^{-1}(u\ast \widetilde{f}\ |_{N})(n,0,a)  \notag \\
&=&u\ast _{c}\widetilde{f}|_{\textfranc }(n,a,0)\in \chi ^{-1}(u\ast 
\widetilde{I}\ |_{N})(n,a,0)  \notag \\
&=&u\ast \widetilde{I}\ |_{\textfranc }(n,a,0)
\end{eqnarray}

\textbf{Corollary 5.1}. \textit{Let }$I$ \textit{be a subspace of the space} 
$L^{1}(N)$ \textit{and} $\tau I$ \textit{its image by the mapping} $\tau $ 
\textit{such that} $J=\tau I|_{\textfranc }$ \textit{is an} \textit{ideal in}
$L^{1}(\textfranc ),$ \textit{then the following conditions are verified}.

$(1)$ $J$ \textit{is a closed ideal in the algebra} $L^{1}(\textfranc )$ 
\textit{if and only if} $I$ \textit{is a left closed ideal in the algebra }$%
L^{1}(N).$

$(2)J$ \textit{is a prime ideal in the algebra} $L^{1}(\textfranc )$ \textit{%
if and only if} $I$ \textit{is a left prime ideal in the algebra }$L^{1}(N)$

$(3)J$ \textit{is a maximal ideal in the algebra} $L^{1}(\textfranc )$ 
\textit{if and only if} $I$ \textit{is a left maximal ideal in the algebra }$%
L^{1}(N)$

$(4)$ $J$ \textit{is a dense ideal in the algebra} $L^{1}(\textfranc )$%
\textit{\ if and only if }$I$\textit{\ is a left dense ideal in the algebra} 
$L^{1}(N).$

The proof of this corollary results immediately from theorem 5.2.

\section{\protect\bigskip \protect\bigskip Left Ideals of the Group Algebra $%
L^{1}(N\rtimes S).$}

\bigskip Let $S=SL(2,\mathbb{R)}/SO(2)$ the symmetric space of the real semi
simple Lie group $SL(2,\mathbb{R)}$, which is diffeomorphism on the group 
\begin{equation}
S=SL(2,\mathbb{R})/SO(2)=\{\left( X=%
\begin{array}{cc}
a & n \\ 
0 & a^{-1}%
\end{array}%
\right) ,a\in \mathbb{R}_{+}^{\ast }\}
\end{equation}

\bigskip The group $S$ is isomorphic onto the group $\mathbb{R}\rtimes
_{\varrho }\mathbb{R}_{+}^{\ast }$ semidirect product of the two group $%
\mathbb{R}$ and $\mathbb{R}_{+}^{\ast }$ where $\varrho :\mathbb{R}%
_{+}^{\ast }\mathbb{\rightarrow }Aut\mathbb{\ (R})$ is the group
homomorphism from the real group into the group $Aut\mathbb{\ (R})$ of all
automorphisms of the vector group $\mathbb{R},$ defined as 
\begin{equation*}
\varrho (x)(n)=xn
\end{equation*}

First, I will prove the solvability of any invariant differential operator
on the connected solvable group $S.$ Therefor denote by $W=\mathbb{R}\times 
\mathbb{R}_{+}^{\ast }\times \mathbb{R}_{+}^{\ast },$ with the following law
defined as%
\begin{equation}
\ (n,x,y)(m,a,b)=(n+\varrho (y)m,xa,yb)=(n+m+y^{2}a,xa,yb)
\end{equation}%
for any $(n,a)\in S$, and $(m,b)\in S$ , here $\varrho (a)m=a^{2}m.$ Let $K$
be the group $\mathbb{R}\times \mathbb{R}_{+}^{\ast }$, which is the direct
product of the group $\mathbb{R}$ with the group $\mathbb{R}_{+}^{\ast }.$
So the group $S$ can be identified with the subgroup $\mathbb{R}\times 
\mathbb{\{}1\mathbb{\}}\times \mathbb{R}_{+}^{\ast }$ of $W$ and $K$ can be
identified with the subgroup $\mathbb{R}\times \mathbb{R}_{+}^{\ast }\times 
\mathbb{\{}1\mathbb{\}}$ of $W.$

\textbf{Definition 6.1}. \textit{For any function} $f\in \mathcal{D}(S),$ 
\textit{we can define a function} $\tau f$ \textit{on }$W$ \textit{by}%
\begin{equation}
\tau f(n,a,b)=f(\varrho (a)n,ab)
\end{equation}

\textbf{Remark 6.1. }The function $\tau f$\ is invariant in the following
sense

\begin{equation}
\tau f(\varrho (x^{-1})n,xa,x^{-1}b)=\tau f(n,a,b)
\end{equation}

\bigskip Now denote by $\tau (C^{\infty }(S))$ $(resp.$ $\tau (C^{\infty
}(K))$ ) the image of $C^{\infty }(S$ $)$ $(resp.C^{\infty }(K)$\ by the
transformation $\tau $, then we have%
\begin{equation*}
\ \tau (C^{\infty }(S))|_{S}=C^{\infty }(S)
\end{equation*}%
\begin{equation}
\ \tau (C^{\infty }(K))|_{K}=C^{\infty }(K)
\end{equation}

\bigskip \textbf{Definition 6.2. }\textit{Let be the mapping} $\chi $ $:\tau
(C^{\infty }(K))|_{K}$ $\longrightarrow \tau (C^{\infty }(S))|_{S}$\textit{\
defined by}

\begin{equation}
\tau f|_{K}\ (z,y,1)\rightarrow \tau f|_{N}(z,1,y)
\end{equation}%
\begin{equation}
\tau f|_{K}\ (n,a,1)\rightarrow \tau f|_{S}(n,1,a)
\end{equation}%
\textit{is topological isomorphisms and its inverse is nothing but }$\chi
^{-1}$\textit{defined by}

\begin{equation}
\tau f|_{S}\ (n,1,a)\rightarrow \tau f|_{K}(n,a,1)
\end{equation}

My main result is

\textbf{Theorem 6.1.} \textit{If }$P_{u}$\ \textit{any} \textit{invariant
differential operator on }$S$\textit{\ associated to the distribution }$u\in 
\mathcal{U}$\textit{, then, we have }

\begin{equation}
P_{u}\text{ }C^{\infty }(S)=C^{\infty }(S)
\end{equation}

\textit{Proof:} Let $Q_{u}$\ be the invariant differential operator with
constant coefficients on $K$ associated to $u$ , then by the theory of
differential operators with constant coefficients $[20]$, we get 
\begin{equation}
Q_{u}\ \ \tau (C^{\infty }(K))|_{K}=\tau (C^{\infty }(K))|_{K}=C^{\infty }(K)
\end{equation}

That means for any $\psi (n,a)\in C^{\infty }(K),$ there exist a function $%
\varphi (n,a,x)\in \tau (C^{\infty }(K))|_{K},$ such that%
\begin{equation}
Q_{u}\varphi (n,a,1)=u\ast _{c}\varphi (n,a,1)=\psi (n,a)
\end{equation}

The function $\psi (n,a)$ can be \ transformed as an invariant function $%
\psi \in \tau (C^{\infty }(K))|_{K}$ as follows

\begin{equation}
\psi (n,a)=\tau \psi (\varrho (a^{-1})n,a,1)
\end{equation}

In other side, we have%
\begin{eqnarray}
&&\chi Q_{u}\ \ \varphi (n,a,1)  \notag \\
&=&Q_{u}\ \ \varphi (n,1,a)=u\ast _{c}\varphi (n,1,a)  \notag \\
&=&u\ast \varphi (n,1,a)=P_{u}\text{ }\varphi (n,1,a)  \notag \\
&=&\chi \tau \psi (\varrho (a^{-1})n,a,1)=\tau \psi (\varrho (a^{-1})n,1,a) 
\notag \\
&=&\psi (n,a)
\end{eqnarray}

So the proof of the solvability of any right invariant differential operator
on $S.$

If $I$ is a subspace of $L^{1}(S),$ we denote by $\tau I$ its image by the
mapping $\tau $, let $\omega =$ $\tau I|_{K}.$ My main result is:

\textbf{Theorem 6.2.} \textit{Let} $I$ \textit{be a subspace of} $L^{1}(S),$ 
\textit{then the following conditions are equivalents}.

$(i)$ $\omega =\tau I\ |_{K}$ \textit{is an ideal in the Banach algebra} $%
L^{1}(K).$

$(ii)$ $I$ \textit{is a left ideal in the Banach algebra} $L^{1}(S).$

\bigskip \textit{Proof:} $(i)$ implies $(ii)$\ Let $I$ be a subspace of the
space $L^{1}(S)$ such that $\omega =\tau I|_{K}$ \ is an ideal in $L^{1}(K),$
then we have: 
\begin{equation}
u\ast _{c}\tau I\ |_{K}(n,a,1)\subseteq \tau I|_{K}(n,a,1)
\end{equation}%
for any $u\in L^{1}(K)$ and $(n,a)\in K$, where%
\begin{equation}
u\ast _{c}\tau I\ |_{K}(n,a,1)=\left\{ \int\limits_{K}\tau f|_{K}\ \left[
n-m,a-b,1)\right] u(m,b)dm\frac{db}{b},\text{ }f\in \text{ }I\right\}
\end{equation}

It shows that%
\begin{equation}
u\ast _{c}\tau f|_{K}(n,a,1)\in \tau I|_{K}(n,a,1)
\end{equation}%
for any $\tau f\in \tau I.$ According to equation$(82),$ we get%
\begin{eqnarray}
&&\chi (u\ast _{c}\tau f|_{K})(n,a,1)  \notag \\
&=&u\text{ }\ast \tau f(n,1,a)\in \chi (\tau I|_{K})(n,a,1)  \notag \\
&=&\tau I\ |_{S}(n,1,a)=I\ (n,a)
\end{eqnarray}%
$(ii)$ implies $(i),$ if $I$ is an ideal in $L^{1}(S),$ then we get 
\begin{eqnarray}
&&u\ast \tau I\ |_{S}(n,1,a)  \notag \\
&=&u\ast I\ (n,a)\subseteq \tau I|_{S}(n,1,a)=I\ (n,a)
\end{eqnarray}%
where%
\begin{equation}
u\ast \tau I|_{S}(n,1,a)=\left\{ \int\limits_{S}\tau f|_{S}\ \left[ \rho
(-b)(n-m),1,a-b\right] u(m,b)dm\frac{db}{b},\text{ }f\in \text{ }I\right\}
\end{equation}

By equation $(36),$we obtain%
\begin{eqnarray}
&&\chi ^{-1}(u\ast \tau f|_{S})(n,1,a)  \notag \\
&=&u\ast _{c}\tau f|_{K}(n,a,1)\in \chi ^{-1}(u\ast \tau I|_{S})(n,a,1) 
\notag \\
&=&u\ast \tau I\ |_{K}(n,a,1)
\end{eqnarray}

\textbf{Corollary 6.1}. \textit{Let }$I$ \textit{be a subspace of the space} 
$L^{1}(S)$ \textit{and} $\tau I$ \textit{its image by the mapping} $\tau $ 
\textit{such that} $\omega =\tau I|_{K}$ \textit{is an} \textit{ideal in} $%
L^{1}(K),$ \textit{then the following conditions are verified}.

$(1)$ $\omega $ \textit{is a closed ideal in the algebra} $L^{1}(K)$ \textit{%
if and only if} $I$ \textit{is a left closed ideal in the algebra }$%
L^{1}(S). $

$(2)$ $\omega $ \textit{is a prime ideal in the algebra} $L^{1}(K)$ \textit{%
if and only if} $I$ \textit{is a left prime ideal in the algebra }$L^{1}(S)$

$(3)$ $\omega $ \textit{is a maximal ideal in the algebra} $L^{1}(K)$ 
\textit{if and only if} $I$ \textit{is a left maximal ideal in the algebra }$%
L^{1}(S)$

$(4)$ $\omega $ \textit{is a dense ideal in the algebra} $L^{1}(K)$\textit{\
if and only if }$I$\textit{\ is a left dense ideal in the algebra} $%
L^{1}(S). $

The proof of this corollary results immediately from theorem \textbf{6.2}.

The Heisenberg group $N$ is the semi-direct product of the two vector Lie
group $\mathbb{R}^{2}\rtimes _{\sigma }\mathbb{R}$. I extend the group $M=N%
\mathbb{\times }S$ by considering the new group $V=$ $\mathbb{R}^{2}\times 
\mathbb{R\times R\times }S$ with the following law%
\begin{eqnarray}
&&X\cdot Y  \notag \\
&=&(n_{3},n_{2},n_{1},n_{4},a_{1},a_{2},a_{3})(m_{3},m_{2},m_{1},m_{4},b_{1},b_{2},b_{3})
\notag \\
&=&((n_{3}+m_{3}+\sigma
(n_{4})(m_{3},m_{2}),n_{2}+m_{2},n_{1}+m_{1},n_{4}+m_{4}),(a_{1}b_{1},a_{2}b_{2},a_{3}b_{3}))
\notag \\
&=&((n_{3}+m_{3}+n_{4}m_{2},n_{2}+m_{2},n_{1}+m_{1},n_{4}+m_{4}),(a_{1}b_{1},a_{2}b_{2},a_{3}b_{3}))
\end{eqnarray}

Denote by $B=$ $\mathbb{R}^{2}\times \mathbb{R\times }S$ the commutative Lie
group of the direct product of three Lie groups $\mathbb{R}^{2},$ $\mathbb{R}
$, and $S.$ In this case the group $M=N\mathbb{\times }S$ can be identified
with the sub-group $\mathbb{R}^{2}\times \mathbb{\{}0\mathbb{\}\times
R\times }S$ and the group $B=\mathbb{R}^{2}\times \mathbb{R\times }S$ can be
identified with the sub-group\ $\mathbb{R}^{2}\mathbb{\times R\times \{}0%
\mathbb{\}\times }S$

\textbf{Definition 6.3. }\textit{Any function }$\psi \in C^{\infty }(M)$ 
\textit{can be extended to a unique function }$\Xi \psi $ \textit{belongs to 
}$C^{\infty }(V),$ as follows%
\begin{eqnarray}
&&\Xi \psi ((n_{3},n_{2},n_{1},n_{4}),s)  \notag \\
&=&\psi ((\sigma (n_{1})(n_{3},n_{2}),n_{1}+n_{4}),s)  \notag \\
&=&\psi ((n_{1}(n_{3},n_{2}),n_{1}+n_{4}),s)  \notag \\
&=&\psi ((n_{3}+n_{1}n_{2},n_{2},n_{1}+n_{4}),s)
\end{eqnarray}%
\textit{for any }$(n_{3},n_{2},n_{1},n_{4})\in N\times \mathbb{R},s\in S,$ $%
n_{1}(n_{3},n_{2})=(n_{3}+n_{1}n_{2},n_{2})=\sigma (n_{1})(n_{3},n_{2})$

If $I$ is a subspace of $L^{1}(M),$ we denote $\Xi I$ its image by the
mapping $\Xi $. Let $J=$ $\Xi I\ |_{B}.$

My main result is:

\textbf{Theorem 6.3.} \textit{Let} $I$ \textit{be a subspace of} $L^{1}(K),$ 
\textit{then the following conditions are equivalents}.

$(i)$ $J=\widetilde{I}\ |_{B}$ \textit{is an ideal in the Banach algebra} $%
L^{1}(B).$

$(ii)$ $I$ \textit{is a left ideal in the Banach algebra} $L^{1}(M).$

For the proof of this theorem, I refer to my book $[9,$ $ChapI,theorem$ $%
3.1.]$\ \ \ 

\textbf{Corollary 6.2}. \textit{Let }$I$ \textit{be a subspace of the space} 
$L^{1}(M)$ \textit{and} $\Xi I$ \textit{its image by the mapping} $\Xi $ 
\textit{such that} $J=\Xi I|_{B}$ \textit{is an} \textit{ideal in} $%
L^{1}(B), $ \textit{then the following conditions are verified}.

$(1)$ $J$ \textit{is an ideal in the algebra} $L^{1}(B)$ \textit{if and only
if} $I$ \textit{is a closed ideal in the algebra }$L^{1}(M)$ \textit{if and
only if} $I$ \textit{is a closed left ideal in the algebra }$L^{1}(N\rtimes
S).$

$(2)J$ \textit{is a prime ideal in the algebra} $L^{1}(B)$ \textit{if and
only if} $I$ \textit{is a prime ideal in the algebra }$L^{1}(M)$ \textit{if
and only if} $I$ \textit{is a prime left ideal in the algebra }$%
L^{1}(N\rtimes S)$

$(3)J$ \textit{is a maximal ideal in the algebra} $L^{1}(B)$ \textit{if and
only if} $I$ \textit{is a maximal ideal in the algebra }$L^{1}(M)$ \textit{%
if and only if} $I$ \textit{is a left maximal ideal in the algebra }$%
L^{1}(N\rtimes S)$

$(4)$ $J$ \textit{is a dense ideal in the algebra} $L^{1}(B)$\textit{\ if
and only if }$I$\textit{\ is a dense ideal in the algebra} $L^{1}(M)$ 
\textit{if and only if }$I$\textit{\ is a left dense ideal in the algebra} $%
L^{1}(N\rtimes S)$

For the proof of this theorem, I refer to \textbf{Theorem 6.2. and Corollary
6.1}.

\end{document}